\documentclass[10pt, article]{amsart}

\usepackage{tikz}
\usetikzlibrary{calc}
\usepackage{ae} 
\usepackage[T1]{fontenc}
\usepackage[cp1250]{inputenc}
\usepackage{amsmath}
\usepackage{amssymb, amsfonts,amscd,verbatim}

\usepackage[normalem]{ulem}
\usepackage{hyperref}
\usepackage{indentfirst}
\usepackage{latexsym}
\input xy
\xyoption{all}

\usepackage{amsmath}    

\theoremstyle{plain}
\newtheorem{Pocz}{Poczatek}[section]
\newtheorem{Proposition}[Pocz]{Proposition}

\newtheorem{Example}[Pocz]{Example}

\theoremstyle{definition}

\theoremstyle{remark}

\newtheorem{Exercise}[Pocz]{Exercise}

\numberwithin{equation}{section}
\title[Anomalous coverings]
{Anomalous coverings}

\author{Jerzy Dydak}

\address{University of Tennessee, Knoxville, TN 37996}
\email{dydak\@@math.utk.edu}


\date{October 6, 2012.}

\begin{document}
\maketitle
\begin{center}
\today
\end{center}

\begin{abstract}

We give examples of anomalous two-fold coverings $p:E\to B$ of connected spaces:\\
a. one where $B$ is simply connected,\\
b. the other of path-connected spaces that has an Evil Twin; a non equivalent covering
$q:E\to B$ with the same image of the fundamental group.

\end{abstract}

\section{Introduction}

There is recent interest in constructing theories of covering maps in various settings (Berestovskii-Plaut 
\cite{BP3} and Brodskiy-Dydak-LaBuz-Mitra \cite{BDLM2}-\cite{BDLM3} for uniform spaces, Fischer-Zastrow \cite{FisZas} and  Brodskiy-Dydak-LaBuz-Mitra \cite{BDLM} for locally path-connected spaces, and Dydak \cite{Dyd} for general spaces). Those efforts amount to generalizing the concept of coverings. Another way to proceed (in order to get a workable theory) is to narrow down covering projections. This was done by R.H.Fox \cite{Fox1}, \cite {Fox2}
who created the concept of \textbf{overlays}.

It is well-known that the classical theory of covers works best for locally semi-simple connected spaces that are locally path-connected.
There is an example of Zeeman \cite[6.6.14 on p.258]{HilWyl} that points out the limits of the classical theory. That example amounts to two non-equivalent coverings of non-locally path-connected
spaces
with the same image of the fundamental groups. Yet this example is not mentioned in current textbooks on topology.
\par

The purpose of this note is to outline an "evolutionary" way of arriving at examples of anomalous coverings. Our example of non-equivalent coverings of path-connected spaces is stronger than that of Zeeman \cite[6.6.14 on p.258]{HilWyl} in the sense that
the total spaces are (naturally) homeomorphic. Also, it captures the essential features of Zeeman's example.
\par
The note is written in the style of a Moore School textbook (guiding students to results via a sequence of definitions and problems). Hopefully it will be used for student presentations in topology classes. Notice there are no proofs included - a promising student ought to be able to reconstruct them on his/her own.
\par
Be aware that we are coining a few new terms. While Sine Curve and Warsaw Circle are widely used, Dusty Broom and Zeeman's Palm seem to be brand new.

The author is grateful to Greg Conner for pointing out the need to cite overlays of R.H.Fox \cite{Fox1}, \cite {Fox2}.

\section{Anomalous coverings}

Throughout this section $X$ is a connected space with two path components, each of them simply connected and locally path-connected.

\begin{Example}
$X$ is the \textbf{Sine Curve}, the closure on the plane of the graph of the function $f(x)=\sin(\frac{\pi}{x})$, $0 < x \leq 1$.
\end{Example}

\begin{Example}
$X$ is the \textbf{Dusty Broom}, the union of the infinite broom (the union of straight arcs joining $(0,0)$ and $(1,\frac{1}{n})$, $n\ge 1$) and a speck of dust in the form of the point $(1,0)$.
\end{Example}
The major difference between the Sine Curve and Dusty Broom is that the Sine Curve is compact.

The simplest way to make $X$ path-connected is to add an arc $A$ joining points $x_0$ and $x_1$ from different path-components of $X$ (with the interior of $A$ disjoint from $X$). We shall refer to such an arc as a \textbf{bridge}
joining $x_0$ and $x_1$.

In case of the Sine Curve we join $(1,0)$ and $(0,-1)$ resulting in the \textbf{Warsaw Circle}. In case of the Dusty Broom we join $(0,0)$ with $(1,0)$ resulting in \textbf{Zeeman's Palm} (infinitely many fingers and one thumb).

\begin{Exercise}
Show $X_1=X\cup A$ is path-connected and simply connected.
\end{Exercise}

\begin{Proposition}\label{FirstCovering}
Consider the space $\tilde X_1$ obtained from $X\times \{0,1\}$ by adding two bridges $A_0$ and $A_1$; one joining $(x_0,0)$
and $(x_1,1)$, the other joining $(x_0,1)$ and $(x_1,0)$. The natural extension $p:\tilde X_1\to X_1:=X\cup A$
of the projection $X\times \{0,1\}\to X$ is a two-fold covering map.
\end{Proposition}

\begin{Exercise}
Show $\tilde X_1$ is a connected space with two path components, each of them simply connected.
\end{Exercise}

We are encountering first anomalies in the theory of coverings: \\
1. A connected cover of a path-connected space may not be path-connected.\\
2. A simply connected space $B$ may admit a connected cover $E\to B$ that is not a homeomorphism. 

What happens if we try to improve $\tilde X_1$ by making it path-connected? Our standard way is to add a bridge $B_0$ joining its two path-components - let's make $B_0$ join $(x_0,0)$ and $(x_1,0)$. Now the image of $\tilde X_1\cup B_0$ is $X_2:=X\cup A\cup B$, where $B$ is another bridge joining $x_0$ and $x_1$.

Obviously, we want to extend $p$ to a covering projection, so we need to add one more bridge $B_1$ joining
$(x_0,1)$ and $(x_1,1)$ resulting in a path-connected space $\tilde X_2$ and a two-fold covering map $p_1:\tilde X_2\to X_2$.

\begin{Exercise}
$\pi_1(X_2)=Z$, $\pi_1(\tilde X_2)=Z$, and the image of $\pi_1(p_1):\pi_1(\tilde X_2)\to \pi_1(X_2)$ is $2\cdot Z$.
\end{Exercise}

Notice there is another covering map $p_2: \tilde X_2\to X_2$ (the Evil Twin of $p_1$) extending the projection $X\times \{0,1\}\to X$. Namely, we exchange the bridges: $A_0$ and $A_1$ are sent homeomorphically onto $B$, and the other two bridges $B_0$, $B_1$ are sent homeomorphically onto $A$. Abstractly speaking, it is the same covering map (we simply change the labelling of the bridges in $\tilde X_2$), so the image of $\pi_1(p_2):\pi_1(\tilde X_2)\to \pi_1(X_2)$ is $2\cdot Z$. Yet

\begin{Proposition}
There is no continuous $f:\tilde X_2\to \tilde X_2$ such that $p_i\circ f= p_j$ for $i\ne j$, $1\leq i,j\leq 2$.
\end{Proposition}

\begin{Exercise}
Zeeman \cite[6.6.14 on p.258]{HilWyl} constructs a pair of two-fold coverings over the wedge $B$ of the unit circle $S^1$ and Zeeman's Palm $ZP$. In the first one the total space $E_1$ is $S^1$ with two copies of $ZP$ attached at $1$ and $-1$, respectively. The projection $p_1:E_1\to B$ is the natural extension of the two-fold covering $z\to z^2$ of $S^1$ over itself. In the second one the total space $E_2$ is $S^1$ with two copies of the Dusty Broom attached at $1$ and $-1$, respectively. One then adds two bridges: each from the base of one broom to the speck of dust of the other. The projection $p_2:E_2\to B$ is the natural extension of the two-fold covering $z\to z^2$ of $S^1$ over itself.
Show $E_1$ and $E_2$ are not homeomorphic.
\end{Exercise}

\section{Reflection}
Let's reflect on why we considered spaces $X$ that are connected with two path components, each of them simply connected and locally path-connected.

\begin{Exercise}
Suppose $p:E\to B$ is a two-fold covering of connected spaces. If $B$ is path-connected and simply connected, show $E$ has exactly two path components, they are homeomorphic, and each of them is simply connected.
\end{Exercise}

The last exercise is to check if the reader understands material.
\begin{Exercise}
Construct a two-fold covering $p:E\to B$ of connected spaces with the following properties:\\
a. $B$ is path-connected and simply connected,\\
b. path-components of $E$ are locally path-connected.
\end{Exercise}


\begin{thebibliography}{99}
\bibitem{BP3}
V. Berestovskii, C. Plaut, {\em Uniform universal covers of uniform spaces},
Topology Appl. 154 (2007), 1748--1777.

 \bibitem{BDLM} N.Brodskiy, J.Dydak, B.Labuz, A.Mitra,
{\em Covering maps for locally path-connected spaces}, Fundamenta Mathematicae 218 (2012), 13--46.
\par\noindent
{\tt\verb+http://front.math.ucdavis.edu/0801.4967+}

\bibitem{BDLM2} N.Brodskiy, J.Dydak, B.Labuz, A.Mitra,
{\em Rips complexes and covers in the uniform category}, 
{\tt\verb+http://front.math.ucdavis.edu/0706.3937+}

   \bibitem{BDLM3} N. Brodskiy, J. Dydak, B. Labuz, A. Mitra, \emph{Topological and uniform structures on universal covering spaces}, arXiv:1206.0071

   \bibitem{Dyd} J. Dydak, \emph{Coverings and fundamental groups: a new approach}, arXiv:1108.3253

\bibitem{FisZas}  H.Fischer, A.Zastrow, {\em Generalized universal
coverings and the shape group}, Fundamenta Mathematicae 197 (2007), 167--196.

\bibitem{Fox1} R. H. Fox, \emph{On shape}, Fund.Math. 74 (1972), 47--71.
\bibitem{Fox2} R. H. Fox, \emph{Shape theory and covering spaces}, Lecture Notes in Math., Vol. 375, Springer,
Berlin, 1974, pp 77--90.

\bibitem{HilWyl}
P.J. Hilton, S. Wylie, {\em Homology theory: An introduction to algebraic topology}, Cambridge University Press, New York 1960 xv+484 pp.
\end{thebibliography}
\end{document}